 \documentclass[12pt]{article}

 \usepackage{mathrsfs,amsfonts,amssymb}
 \usepackage[dvips]{color}

 \newtheorem{lemma}{Lemma}[section]
 \newtheorem{proposition}{Proposition}[section]
 \newtheorem{theorem}{Theorem}[section]
 \newtheorem{corollary}{Corollary}[section]

 \def\blemma{\begin{lemma}\sl{}\def\elemma{\end{lemma}}}
 \def\bproposition{\begin{proposition}\sl{}\def\eproposition{\end{proposition}}}
 \def\btheorem{\begin{theorem}\sl{}\def\etheorem{\end{theorem}}}
 \def\bcorollary{\begin{corollary}\sl{}\def\ecorollary{\end{corollary}}}

 \def\beqlb{\begin{eqnarray}}\def\eeqlb{\end{eqnarray}}
 \def\beqnn{\begin{eqnarray*}}\def\eeqnn{\end{eqnarray*}}

 \def\ar{\!\!\!&}\def\nnm{\nonumber}

 \def\<{\langle}\def\>{\rangle}

 \def\mbb{\mathbb}\def\mbf{\mathbf}\def\mcr{\mathscr}

 \def\proof{\noindent\textit{Proof.~~}}
 \def\qed{\hfill$\square$\medskip}

\begin{document}

\noindent{Published in: \textit{Stochastic Analysis and
Applications} \textbf{25} (2007), 6: 1273--1296}

\bigskip\bigskip

\centerline{\LARGE\bf Continuous local time of a purely}

\medskip

\centerline{\LARGE\bf atomic immigration superprocess}

\medskip
\centerline{\LARGE\bf with dependent spatial motion}

\bigskip\bigskip

\centerline{\bf Zenghu Li\,$^1$ and Jie Xiong\,$^{2,3}$}

\bigskip

\centerline{\small $^1$ School of Mathematical Sciences, Beijing
Normal University,}

\centerline{\small Beijing 100875, P.R. China. E-mail: \tt
lizh@bnu.edu.cn}

\bigskip

\centerline{\small $^2$ Department of Mathematics, University of
Tennessee, Knoxville,}

\centerline{\small TN 37996-1300, U.S.A. E-mail: \tt
jxiong@math.utk.edu}

\bigskip

\centerline{\small $^3$ Department of Mathematics, Hebei Normal
University,}

\centerline{\small Shijiazhuang 050016, P.R. China}

\bigskip\bigskip

{\narrower{\narrower

\centerline{\large\bf Abstract}

\bigskip

A purely atomic immigration superprocess with dependent spatial
motion in the space of tempered measures is constructed as the
unique strong solution of a stochastic integral equation driven by
Poisson processes based on the excursion law of a Feller branching
diffusion, which generalizes the work of Dawson and Li \cite{DL03}.
As an application of the stochastic equation, it is proved that the
superprocess possesses a local time which is H\"older continuous of
order $\alpha$ for every $\alpha < 1/2$. We establish two scaling
limit theorems for the immigration superprocess, from which we
derive scaling limits for the corresponding local time.

\bigskip

\textit{Key words}: local time, superprocess, dependent spatial
motion, immigration, excursion, Poisson random measure; scaling
limit theorem.

\bigskip

\textit{Mathematics Subject Classification (2000)}: 60J80; 60G57;
60H20

\bigskip

\textit{Running head}: Continuous local time of immigration
superprocess.

\par}\par}

\bigskip\bigskip


\section{Introduction}

\setcounter{equation}{0}

Let $p\ge0$ and let $\phi_p(x) = (1+x^2)^{-p/2}$ for $x\in \mbb{R}$.
We denote by $C_p(\mbb{R})$ the set of continuous functions $\phi$
on $\mbb{R}$ satisfying $|\phi| \le$ const$\,\cdot\,\phi_p$ and
denote by $C_p^2(\mbb{R})$ the subset of $C_p(\mbb{R})$ consisting
of twice continuously differentiable functions $\phi$ with
$|\phi^\prime|+|\phi^{\prime\prime}| \le$ const$\,\cdot\,\phi_p$.
Let $M_p(\mbb{R})$ denote the space of tempered Borel measures $\mu$
on $\mbb{R}$ such that
 \beqnn
\<\phi,\mu\> := \int_{\mbb{R}}\phi(x) \mu(dx) <\infty
 \eeqnn
for every $\phi\in C_p(\mbb{R})$. Let $M_p^a(\mbb{R})$ be the subset
of $M_p(\mbb{R})$ consisting of purely atomic measures. In the case
$p=0$, we simply write $C(\mbb{R})$ and $M(\mbb{R})$ instead of
$C_0(\mbb{R})$ and $M_0(\mbb{R})$, respectively. Let ``$\|\cdot\|$''
denote the supremum norm. Suppose that $h$ is a continuously
differentiable function on $\mbb{R}$ such that both $h$ and
$h^\prime$ are square-integrable. Then the function
 \beqlb\label{1.1}
\rho(x) = \int_{\mbb{R}}h(y-x)h(y) dy, \qquad x\in\mbb{R}
 \eeqlb
is twice continuously differentiable with bounded derivatives
$\rho^\prime$ and $\rho^{\prime\prime}$. We fix a constant
$\sigma>0$ and a measure $m\in M_p(\mbb{R})$. Let $q(\nu,a)$ be a
Borel function of $(\nu,a) \in M_p(\mbb{R}) \times \mbb{R}$
satisfying certain regularity conditions to be specified. A
martingale problem for a continuous process $\{Y_t: t\ge0\}$ in
$M_p(\mbb{R})$ can be formulated in the following way: For each
$\phi \in C_p^2(\mbb{R})$,
 \beqlb\label{1.2}
M_t(\phi) = \<\phi,Y_t\> - \<\phi,Y_0\> - \frac{1}{\,2\,}\rho(0)
\int_0^t \<\phi^{\prime\prime},Y_s\> ds - \int_0^t \<q(Y_s,\cdot)
\phi, m\> ds
 \eeqlb
is a continuous martingale with quadratic variation process
 \beqlb\label{1.3}
\<M(\phi)\>_t = \int_0^t\<\sigma\phi^2,Y_s\> ds + \int_0^t ds
\int_{\mbb{R}^2} \rho(x-y)\phi^\prime(x)\phi^\prime(y)
Y_s(dx)Y_s(dy).
 \eeqlb
A solution $\{Y_t: t\ge0\}$ of the martingale problem is called an
\textit{immigration superprocess with dependent spatial motion}
(ISDSM); see Ref.~\cite{DL03}. The phrase \textit{superprocess with
dependent spatial motion} (SDSM) naturally refers to the special
case $m=0$. Compared with the classical super Brownian motion over
the real line, the last term in (\ref{1.2}) represents an
immigration factor with interactive immigration rate while the last
term in (\ref{1.3}) comes from the dependent spatial motion. In
particular, if $q(\nu,a) = q(a)$ only depends on $a\in \mbb{R}$, the
immigration becomes non-interactive and the uniqueness of solution
of the martingale problem can be proved by a duality argument or a
conditional log-Laplace functional; see
Refs.~\cite{DLW01,LLW04a,LWX05,W98}. The uniqueness of solution of
the general martingale problem still remains a challenging open
problem.

An ISDSM was constructed in Ref.~\cite{DL03} in the special case
$p=0$. Instead of the martingale problem, the authors considered a
stochastic integral equation driven by Poisson random measures based
on the excursion law of a Feller branching diffusion. They showed
that there is a unique strong solution of the equation which also
solves (\ref{1.2}) and (\ref{1.3}); see also Refs.~\cite{FL04,S90}.
In this paper we extend the result of Ref.~\cite{DL03} to the
general state space $M_p(\mbb{R})$. For the study of the ISDSM, the
stochastic equation has several advantages over the martingale
problem formulation. For instance, the stochastic equation provides
much more information on the structures of the ISDSM than the
martingale problem. As pointed out in Ref.~\cite{DL03}, from the
equation we know immediately that the ISDSM lives in the space of
purely atomic measures in contrast to the classical super Brownian
motion; see also Ref.~\cite{W97}. With that observation, it is
natural to ask whether or not the corresponding occupation time
process
 \beqlb
Z_t := \int_0^t Y_sds, \qquad t\ge0
 \eeqlb
is absolutely continuous with respect to the Lebesgue measure. We
give an answer to this question as an application of the stochastic
equation. We shall see that $\{Z_t: t\ge0\}$ is really absolutely
continuous and its density field $\{z(b,t): b\in \mbb{R}, t\ge0\}$
can be represented in terms of stochastic integrals of the
excursions with respect to the Poisson random measures and certain
Brownian local times. By this representation we prove that
$\{z(b,t): b\in \mbb{R}, t\ge0\}$ is H\"older continuous of order
$\alpha$ for every $\alpha < 1/2$. It seems difficult to establish
those results from (\ref{1.2}) and (\ref{1.3}). The simple
derivations of the results given here show the efficiency of the
stochastic equation in the study of properties of the immigration
superprocess. As another application of the stochastic equation, we
prove two scaling limit theorems for $\{Y_t: t\ge0\}$. {From} those
theorems we derive the limit theorems for the density field
$\{z(b,t): b\in \mbb{R}, t\ge0\}$.


\section{Preliminary results}

\setcounter{equation}{0}

Let $m \in M_p(\mbb{R})$ and let $(\Omega, \mcr{F}, \mcr{F}_t,
\mbf{P})$ be a filtered standard probability space satisfying the
usual hypotheses. A function $\eta(\cdot,\cdot,\cdot)$ on
$[0,\infty)\times \mbb{R}\times \Omega$ is said to be
\textit{simple} if it is of the form
 \beqlb\label{2.1}
\eta(s,x,\omega) = \eta_0(x,\omega)1_{[r_0,r_1]}(s) +
\sum_{i=1}^\infty \eta_i(x,\omega)1_{(r_i,r_{i+1}]}(s),
 \eeqlb
where $0=r_0<r_1<r_2<\ldots$ and $\eta_i(\cdot,\cdot)$ is
$\mcr{B}(\mbb{R})\times \mcr{F}_{r_i}$-measurable. Let $\mcr{P}$
be the completion with respect to $dsm(da) \mbf{P}(d\omega)$ of
the $\sigma$-algebra on $[0,\infty)\times \mbb{R}\times \Omega$
generated by all simple functions. We say a function on
$[0,\infty)\times \mbb{R}\times \Omega$ is \textit{predictable} if
it is $\mcr{P}$-measurable. Let $L^2_+(\mcr{P})$ denote the set of
all non-negative predictable functions $\eta(\cdot,\cdot,\cdot)$
on $[0,\infty)\times \mbb{R}\times \Omega$ such that
 \beqlb\label{2.2}
\int_0^t\mbf{E}[\<\eta(s,\cdot)\phi_p,m\>^2]ds < \infty, \qquad
t\ge0.
 \eeqlb
Suppose that $\sigma>0$ and $\eta(\cdot,\cdot,\cdot)\in
L^2_+(\mcr{P})$. Given $\mu \in M_p(\mbb{R})$, we consider the
following martingale problem of a continuous process $\{Y_t:
t\ge0\}$ in $M_p(\mbb{R})$: For each $\phi\in C_p^2(\mbb{R})$,
 \beqlb\label{2.3}
M_t(\phi) = \<\phi, Y_t\> - \<\phi, \mu\> - \frac{1}{2}
\rho(0)\int_0^t \<\phi^{\prime\prime}, Y_s\> ds - \int_0^t
\<\eta(s,\cdot)\phi,m\>ds
 \eeqlb
is a continuous martingale relative to $(\mcr{F}_t)_{t\ge0}$ with
quadratic variation process
 \beqlb\label{2.4}
\<M(\phi)\>_t = \int_0^t\<\sigma\phi^2, Y_s\>ds + \int_0^t
ds\int_{\mbb{R}^2} \rho(x-y) \phi^\prime(x)\phi^\prime(y)
Y_s(dx)Y_s(dy).
 \eeqlb
A solution of (\ref{2.3}) and (\ref{2.4}) can be regarded as an
generalized immigration superprocess with immigration rate given
by the two parameter process $\{\eta(s,a): s\ge0, a\in \mbb{R}\}$.

\bproposition\label{p2.1} Suppose that $\{Y_t: t\ge0\}$ is a
solution of the martingale problem given by (\ref{2.3}) and
(\ref{2.4}). Let $c\ge0$ be a constant such that $|\phi_p^\prime|
+ |\phi_p^{\prime\prime}| \le c\phi_p$ and let $C(t) =
c^2\|\rho\|(16+\|\rho\|t)$. Then we have
 \beqlb\label{2.5}
\mbf{E}[\<\phi_p,Y_t\>]
 \le
G_t(\phi_p) + \frac{1}{2} c\|\rho\|\int_0^t
G_s(\phi_p)\exp\Big\{\frac{1}{2} c\|\rho\|(t-s)\Big\}ds
 \eeqlb
and
 \beqlb\label{2.6}
\mbf{E}\Big[\sup_{0\le s\le t}\<\phi_p,Y_s\>^2\Big]
 \le
H_t(\phi_p) + C(t)\int_0^t H_s(\phi_p)\exp\{C(t)(t-s)\}ds,
 \eeqlb
where
 \beqnn
G_t(\phi_p)
 =
\<\phi_p, \mu\> + \int_0^t\mbf{E}[\<\eta(s,\cdot)\phi_p, m\>]ds
 \eeqnn
and
 \beqnn
H_t(\phi_p)
 =
4\<\phi_p, \mu\>^2 +  4t\int_0^t\mbf{E}[\<\eta(s,\cdot)\phi_p,
m\>^2]ds + 16\sigma\int_0^t\mbf{E}[\<\phi_p, Y_s\>]ds.
 \eeqnn
\eproposition

\proof The first inequality follows by taking the expectations in
(\ref{2.3}) and applying Gronwall's inequality. By (\ref{2.3}) and
the martingale inequality we have
 \beqnn
\mbf{E}\Big[\sup_{0\le s\le t}\<\phi_p,Y_s\>^2\Big]
 \ar\le\ar
4\<\phi_p, \mu\>^2 + 4\mbf{E}\bigg[\bigg(\int_0^t
\<\eta(s,\cdot)\phi_p, m\>ds\bigg)^2\bigg] \\
 \ar \ar
+\, \|\rho\|^2 \mbf{E}\bigg[\bigg(\int_0^t
\<\phi_p^{\prime\prime}, Y_s\>ds\bigg)^2\bigg] \\
 \ar \ar
+\, 16\int_0^t\mbf{E}\Big[\sigma\<\phi_p^2, Y_s\> +
\|\rho\|\<\phi_p^\prime, Y_s\>^2\Big]ds  \\
 \ar\le\ar
4\<\phi_p, \mu\>^2 + 4t\int_0^t\mbf{E}[\<\eta(s,\cdot)\phi_p,
m\>^2]ds \\
 \ar \ar
+\, c^2\|\rho\|^2t\int_0^t\mbf{E}[\<\phi_p, Y_s\>^2]ds \\
 \ar \ar
+\, 16\int_0^t\Big(\sigma\mbf{E}[\<\phi_p, Y_s\>] +
c^2\|\rho\|\mbf{E}[\<\phi_p, Y_s\>^2]\Big)ds.
 \eeqnn
Therefore, we can use Gronwall's inequality again to get
(\ref{2.6}). \qed

Clearly, the family of martingales $\{M_t(\phi)\}$ in (\ref{2.3})
defines a martingale measure $M(ds,dx)$ in the sense of
Ref.~\cite{W86}. The next result follows by standard arguments.

\bproposition\label{p2.2} Suppose that $\{Y_t: t\ge0\}$ is a
solution of the martingale problem given by (\ref{2.3}) and
(\ref{2.4}). Then for any $t\ge0$ and $\phi \in C^2_p(\mbb{R})$ we
have a.s.\
 \beqlb\label{2.7}
\<\phi, Y_t\> = \<P_t\phi, \mu\> + \int_0^t \<\eta(s,\cdot)
P_{t-s}\phi, m\>ds + \int_0^t\int_{\mbb{R}} P_{t-s}\phi(x)
M(ds,dx),
 \eeqlb
where $(P_t)_{t\ge0}$ is the semigroup of the Brownian motion with
quadratic variation $\rho(0)dt$. Consequently,
 \beqlb\label{2.8}
\mbf{E}[\<\phi, Y_t\>] = \<P_t\phi, \mu\> + \int_0^t
\mbf{E}[\<\eta(s,\cdot) P_{t-s}\phi, m\>]ds.
 \eeqlb
\eproposition

Now we can give the construction of a solution of the martingale
problem (\ref{2.3}) and (\ref{2.4}) with initial state $\mu \in
M_p^a(\mbb{R})$. Let $\{B(t): t\ge0\}$ be a standard Brownian
motion. For any initial condition $\xi(0)=x\ge0$ the stochastic
differential equation
 \beqlb\label{2.9}
d \xi(t) = \sqrt{\sigma\xi(t)} d B(t), \qquad t\ge0,
 \eeqlb
has a unique solution, which is known as a Feller branching
diffusion with constant branching rate $\sigma$. In the sequel, we
simply call $\{\xi(t): t\ge0\}$ a \textit{$\sigma$-branching
diffusion}. Let $W = C([0,\infty), \mbb{R}^+)$ and let $\tau_0(w) =
\inf\{s>0: w(s)=0\}$ for $w\in W$. We denote by $W_0$ be the set of
paths $w \in W$ such that $w(t)=w(0)=0$ for every $t\ge\tau_0(w)$.
Let $(\mcr{B}(W_0), \mcr{B}_t(W_0))$ be the natural
$\sigma$-algebras on $W_0$ generated by the coordinate process and
let $\mbf{Q}_\kappa$ be the excursion law of the $\sigma$-branching
diffusion defined in Ref.~\cite{DL03}. Suppose that on a complete
standard probability space $(\Omega,\mcr{F},\mbf{P})$ the following
are defined:
\begin{itemize}

\item[(2.a)]
a white noise $W(ds,dy)$ on $[0,\infty)\times \mbb{R}$ based on
the Lebesgue measure;

\item[(2.b)]
a sequence of independent $\sigma$-branching diffusions
$\{\xi_i(t): t\ge0\}$ with deterministic initial values $\xi_i(0)
= \xi_i$, $i=1,2,\cdots$;

\item[(2.c)]
a Poisson random measure $N_1(ds,da,du,dw)$ on $[0,\infty) \times
\mbb{R}\times [0,\infty) \times W_0$ with intensity $dsm(da)du
\mbf{Q}_\kappa (dw)$.

\end{itemize}
We assume that $\{W(ds,dy)\}$, $\{\xi_i(t)\}$ and
$\{N_1(ds,da,du,dw)\}$ are independent of each other. For $t\ge0$
let $\mcr{F}_t$ be the $\sigma$-algebra generated by all
$\mbf{P}$-null sets and the families of random variables
 \beqlb\label{2.10}
\{W([0,s]\times B), \xi_i(s): 0\le s\le t; B\in\mcr{B}(\mbb{R}),
i=1,2,\cdots\}
 \eeqlb
and
 \beqlb\label{2.11}
\{N_1(J\times A): J\in\mcr{B}([0,s]\times \mbb{R}\times [0,\infty));
A\in\mcr{B}_{t-s}(W_0); 0\le s\le t\}.
 \eeqlb
It is known that for any $(r,a) \in [0,\infty) \times \mbb{R}$
there is a unique solution $\{x_{r,a}(t): t\ge r\}$ of the
equation
 \beqlb\label{2.12}
x(t) = a + \int_r^t\int_{\mbb{R}} h(y-x(s))W(ds,dy), \qquad t\ge
r;
 \eeqlb
see Refs.~\cite{DLW01,W97}. Indeed, each $\{x_{r,a}(t): t\ge r\}$ is
a Brownian motion with quadratic variation $\rho(0)dt$. Let $\{a_i:
i=1,2,\cdots\} \subset \mbb{R}$ and assume $\sum_{i=1}^\infty
\xi_i\phi_p(a_i)< \infty$. For $\eta(\cdot,\cdot,\cdot)\in
L^2_+(\mcr{P})$ we define the purely atomic measure-valued process
$\{Y_t: t\ge0\}$ by
 \beqlb\label{2.13}
Y_t = \sum_{i=1}^\infty \xi_i(t)\delta_{x_{0,a_i}(t)} +
\int_0^t\int_{\mbb{R}}\int_0^{\eta(s,a)}\int_{W_0} w(t-s)
\delta_{x_{s,a}(t)}N_1(ds,da,du,dw).
 \eeqlb
(We make the convention that $\int_a^b = \int_{(a,b]}$.)

\btheorem\label{t2.1} The process $\{Y_t: t\ge0\}$ defined by
(\ref{2.13}) has a continuous modification in $M_p^a(\mbb{R})$
which solves the martingale problem given by (\ref{2.3}) and
(\ref{2.4}). \etheorem

\proof For each integer $k\ge1$ let $\eta_k(t,a) = \eta(t,a)\land
(k\phi_p(a))$. Then $t\mapsto \mbf{E}[\<\eta_k(t,\cdot), m\>^2]$
is a locally bounded function and
 \beqnn
\int_0^t\mbf{E}[\<|\eta(s,\cdot)-\eta_k(s,\cdot)|\phi_p,m\>^2]ds
\to 0
 \eeqnn
as $k\to \infty$ for every $t\ge 0$. By Theorem~5.1 of
Ref.~\cite{DL03}, the process $\{Y_t^{(k)}: t\ge0\}$ defined by
 \beqnn
Y_t^{(k)} = \sum_{|a_i|\le k} \xi_i(t)\delta_{x_{0,a_i}(t)} +
\int_0^t\int_{\mbb{R}}\int_0^{\eta_k(s,a)}\int_{W_0} w(t-s)
\delta_{x_{s,a}(t)}N_1(ds,da,du,dw)
 \eeqnn
has a continuous modification in $M_0^a(\mbb{R})$ and hence in
$M_p^a(\mbb{R})$. Note also that $\{Y_t^{(k)}: t\ge0\}$ is
increasing in $k\ge1$. For $n\ge k\ge 1$ let $g_{k,n} = \eta_n
-\eta_k$, $\mu_{k,n} = \sum_{k<|a_i|\le n} \xi_i(0)\delta_{a_i}$
and
 \beqnn
X_t^{(k,n)} = \sum_{k<|a_i|\le n} \xi_i(t)\delta_{x_{0,a_i}(t)} +
\int_0^t\int_{\mbb{R}}\int_{\eta_k(s,a)}^{\eta_n(s,a)} \int_{W_0}
w(t-s) \delta_{x_{s,a}(t)}N_1(ds,da,du,dw).
 \eeqnn
As in the proofs of Theorems~4.1 and~4.5 in Ref.~\cite{DL03} it
follows that $\{X_t^{(k,n)}: t\ge0\}$ has a continuous modification
in $M_p^a(\mbb{R})$ and for each $\phi\in C_p^2(\mbb{R})$,
 \beqnn
M_t^{(k,n)}(\phi) = \<\phi, X_t^{(k,n)}\> - \<\phi, \mu_{k,n}\> -
\int_0^t \<g_{k,n}(s,\cdot)\phi,m\>ds - \frac{1}{2} \rho(0)\int_0^t
\<\phi^{\prime\prime}, X_s^{(k,n)}\> ds
 \eeqnn
is a continuous martingale with quadratic variation process
 \beqnn
\int_0^t\<\sigma\phi^2, X_s^{(k,n)}\>ds + \int_0^t
ds\int_{\mbb{R}^2} \rho(x-y) \phi^\prime(x)\phi^\prime(y)
X_s^{(k,n)}(dx)X_s^{(k,n)}(dy);
 \eeqnn
see also Lemma~3.1 of Ref.~\cite{FL04}. Let $C(t)$ be defined as in
Proposition~\ref{p2.1}. We have
 \beqnn
\mbf{E}[\<\phi_p,X_t^{(k,n)}\>]
 \le
G_t^{(k,n)}(\phi_p) + \frac{1}{2} c\|\rho\|\int_0^t
G_s^{(k,n)}(\phi_p)\exp\Big\{\frac{1}{2} c\|\rho\|(t-s)\Big\}ds
 \eeqnn
and
 \beqnn
\mbf{E}\Big[\sup_{0\le s\le t}\<\phi_p,X_s^{(k,n)}\>^2\Big]
 \le
H_t^{(k,n)}(\phi) + C(t)\int_0^t
H_s^{(k,n)}(\phi)\exp\{C(t)(t-s)\}ds,
 \eeqnn
where
 \beqnn
G_t^{(k,n)}(\phi_p)
 =
\<\phi_p, \mu_{k,n}\> + \int_0^t\mbf{E}[\<g_{k,n}(s,\cdot)\phi_p,
m\>]ds
 \eeqnn
and
 \beqnn
H_t^{(k,n)}(\phi_p)
 \ar=\ar
4\<\phi_p, \mu_{k,n}\>^2 + 4t\int_0^t\mbf{E}[\<g_{k,n}(s,\cdot)
\phi_p, m\>^2]ds   \\
 \ar \ar
+\, 16\sigma\int_0^t\mbf{E}[\<\phi_p, X_s^{(k,n)}\>]ds.
 \eeqnn
Clearly, we have $G_t^{(k,n)}(\phi_p) \to 0$ and $H_t^{(k,n)}
(\phi_p) \to 0$ as $k,n\to\infty$. It follows that
 \beqnn
\mbf{E}\Big[\sup_{0\le s\le t}\<\phi_p,X_s^{(k,n)}\>^2\Big] \to 0
 \eeqnn
as $k,n\to \infty$. Then there is a continuous process $\{X_t:
t\ge0\}$ in $M_p^a(\mbb{R})$ such that
 \beqnn
\mbf{E}\Big[\sup_{0\le s\le t}\<\phi_p,X_s - Y_s^{(k)}\>^2\Big] \to
0
 \eeqnn
as $n\to \infty$ for every $t\ge0$. Clearly, $\{X_t: t\ge0\}$ is
independent of the particular choice of the approximating sequence
$\{\eta_k\}\subset L^2_+(\mcr{P})$ and the martingale
characterizations (\ref{2.3}) and (\ref{2.4}) hold with $\{Y_t:
t\ge0\}$ replaced by $\{X_t: t\ge0\}$.  Then, to finish the proof we
only need to show a.s.\ $X_t = Y_t$ for every $t\ge0$. Recall that
$\int_0^{\eta(s,a)}$ means $\int_{(0,\eta(s,a)]}$ in (\ref{2.13}).
Then we have a.s.\ $X_t \le Y_t$ from the limit procedure for the
construction of $X_t$. Let $\{X_t^{(k)}: t\ge0\}$ be constructed by
the same procedure for $\eta(\cdot,\cdot,\cdot) + 1/k \in
L^2_+(\mcr{P})$. Obviously, we have a.s.\ $Y_t \le X_t^{(k)}$.
However, by Proposition~\ref{p2.2} we have
 \beqnn
\mbf{E}[\<\phi, X_t\>] = \<P_t\phi, \mu\> + \int_0^t
\mbf{E}[\<\eta(s,\cdot)P_{t-s}\phi, m\>]ds
 \eeqnn
and
 \beqnn
\mbf{E}[\<\phi, X_t^{(k)}\>] = \<P_t\phi, \mu\> + \int_0^t
\mbf{E}[\<(\eta(s,\cdot)+1/k)P_{t-s}\phi, m\>]ds.
 \eeqnn
It follows that
 \beqnn
\mbf{E}[\<\phi, X_t\>] \le \mbf{E}[\<\phi, Y_t\>] \le \lim_{k\to
\infty} \mbf{E}[\<\phi, X_t^{(k)}\>] = \mbf{E}[\<\phi, X_t\>].
 \eeqnn
That proves $X_t = Y_t$ a.s.\ for every $t\ge0$. \qed

A solution of (\ref{2.3}) and (\ref{2.4}) with general initial state
$\mu \in M_p(\mbb{R})$ can be constructed in the following way. Let
$\{W(ds,dy)\}$ and $\{N_1(ds,da,du,dw)\}$ be given as in (2.a) and
(2.c) and suppose we are also given
\begin{itemize}

\item[(2.d)]
a Poisson random measure $N_0(da,dw)$ on $\mbb{R}\times W_0$ with
intensity $\mu(da)\mbf{Q}_\kappa (dw)$.

\end{itemize}
We assume that $\{W(ds,dy)\}$, $\{N_0(da,dw)\}$ and
$\{N_1(ds,da,du,dw)\}$ are independent of each other. For $t\ge0$
let $\mcr{F}_t$ be the $\sigma$-algebra generated by all
$\mbf{P}$-null sets and the families of random variables
 \beqlb\label{2.14}
\{W([0,s]\times B), N_0(F\times A): 0\le s\le t, B\in
\mcr{B}(\mbb{R}), A\in\mcr{B}_t(W_0)\}
 \eeqlb
and
 \beqlb\label{2.15}
\{N_1(J\times A): J\in\mcr{B}([0,s]\times \mbb{R}\times [0,\infty)),
A\in\mcr{B}_{t-s}(W_0), 0\le s\le t\}.
 \eeqlb
Let $\{X_t: t\ge 0\}$ be defined by $X_0=\mu$ and
 \beqlb\label{2.16}
X_t = \int_{\mbb{R}}\int_{W_0} w(t) \delta_{x_{0,a}(t)}
N_0(da,dw), \qquad t> 0.
 \eeqlb
For $\eta(\cdot,\cdot,\cdot)\in L^2_+ (\mcr{P})$ we define
 \beqlb\label{2.17}
Y_t = X_t + \int_0^t\int_{\mbb{R}}\int_0^{\eta(s,a)}\int_{W_0}
w(t-s) \delta_{x_{s,a}(t)}N_1(ds,da,du,dw), \qquad t\ge 0.
 \eeqlb
By similar arguments as in the proof of Theorem~\ref{t2.1} we
obtain

\btheorem\label{t2.2} The process $\{Y_t: t\ge 0\}$ defined by
(\ref{2.17}) has a continuous modification in $M_p(\mbb{R})$ which
solves the martingale problem given by (\ref{2.3}) and
(\ref{2.4}). \etheorem


\section{Stochastic equations for superprocesses}

\setcounter{equation}{0}

In this section, we give the construction of the ISDSM by solving
stochastic equations driven by Poisson random measures on the
space of excursions. For any $\mu$ and $\nu \in M_p(\mbb{R})$ set
 \beqlb\label{3.1}
\|\mu-\nu\|_p = \sup_{f\in B_1(\mbb{R})} \bigg|\int_{\mbb{R}}
f(x)\phi_p(x) \mu(dx) - \int_{\mbb{R}} f(x)\phi_p(x)
\nu(dx)\bigg|,
 \eeqlb
where $B_1(\mbb{R})$ denotes the set of Borel functions $f$ on
$\mbb{R}$ such that $|f(x)|\le 1$ for all $x\in \mbb{R}$. Suppose
that $q(\cdot,\cdot)$ is a Borel function on $M_p(\mbb{R})\times
\mbb{R}$ such that there is a constant $K$ such that
 \beqlb\label{3.2}
\<q(\mu,\cdot)\phi_p,m\> \le K(1 + \<\phi_p,\mu\>), \qquad \mu\in
M(\mbb{R}),
 \eeqlb
and for each $R>0$ there is a constant $L_R>0$ such that
 \beqlb\label{3.3}
\<|q(\mu,\cdot) - q(\nu,\cdot)|\phi_p, m\> \le L_R \|\mu-\nu\|_p
 \eeqlb
for $\mu$ and $\nu\in M_p(\mbb{R})$ satisfying $\<\phi_p,\mu\> \le
R$ and $\<\phi_p,\nu\>\le R$. Given $\mu \in M_p(\mbb{R})$, we
consider the following martingale problem of a continuous process
$\{Y_t: t\ge0\}$ in $M_p(\mbb{R})$: For each $\phi\in
C_p^2(\mbb{R})$,
 \beqlb\label{3.4}
M_t(\phi) = \<\phi, Y_t\> - \<\phi, \mu\> - \frac{1}{2}
\rho(0)\int_0^t \<\phi^{\prime\prime}, Y_s\> ds - \int_0^t
\<q(Y_s,\cdot)\phi,m\>ds
 \eeqlb
is a continuous martingale with quadratic variation process
 \beqlb\label{3.5}
\<M(\phi)\>_t = \int_0^t\<\sigma\phi^2, Y_s\>ds + \int_0^t
ds\int_{\mbb{R}^2} \rho(x-y) \phi^\prime(x)\phi^\prime(y)
Y_s(dx)Y_s(dy).
 \eeqlb
By similar arguments as in the proof of Proposition~\ref{p2.1} we
have the following

\bproposition\label{p3.1} Suppose that $\{Y_t: t\ge0\}$ is a
solution of the martingale problem given by (\ref{3.4}) and
(\ref{3.5}). Let $c\ge0$ be a constant such that $|\phi_p^\prime|
+ |\phi_p^{\prime\prime}| \le c\phi_p$ and let $C_1 = K + c
\|\rho\|/2$ and $C_2(t) = c^2\|\rho\| (16+\|\rho\|t)$. Then we
have
 \beqlb\label{3.6}
\mbf{E}[\<\phi_p,Y_t\>]
 \le
G_t(\phi_p) + C_1\int_0^t G_s(\phi_p)\exp\big\{C_1(t-s)\big\}ds
 \eeqlb
and
 \beqlb\label{3.7}
\mbf{E}\Big[\sup_{0\le s\le t}\<\phi_p,Y_s\>^2\Big]
 \le
H_t(\phi_p) + C_2(t)\int_0^t H_s(\phi_p)\exp\{C_2(t)(t-s)\}ds,
 \eeqlb
where $G_t(\phi_p) = \<\phi_p, \mu\> + Kt$ and
 \beqnn
H_t(\phi_p)
 =
4\<\phi_p, \mu\>^2 + 4K^2t^2 + 8(Kt + 2\sigma)
\int_0^t\mbf{E}[\<\phi_p, Y_s\>]ds.
 \eeqnn
\eproposition

Let $\{W(ds,dy)\}$, $\{\xi_i(t)\}$ and $\{N_1(ds,da,du,dw)\}$ be
given as in (2.a), (2.b) and (2.c). Let $(\mcr{F}_t)_{t\ge0}$ be the
filtration generated by (\ref{2.10}) and (\ref{2.11}). By a
\textit{solution} of the stochastic equation
 \beqlb\label{3.8}
Y_t = \sum_{i=1}^\infty \xi_i(t)\delta_{x_{0,a_i}(t)} +
\int_0^t\int_{\mbb{R}}\int_0^{q(Y_s,a)}\int_{W_0} w(t-s)
\delta_{x_{s,a}(t)}N_1(ds,da,du,dw),
 \eeqlb
we mean an $(\mcr{F}_t)$-adapted continuous process $\{Y_t: t\ge
0\}$ in $M_p^a(\mbb{R})$ that satisfies (\ref{3.8}) for all $t\ge
0$.

\btheorem\label{t3.1} There is an unique solution $\{Y_t: t\ge 0\}$
of (\ref{3.8}), which is a strong Markov process in $M_p^a(\mbb{R})$
and solves the martingale problem given by (\ref{3.4}) and
(\ref{3.5}). \etheorem

\proof Suppose that $\{Y_t: t\ge 0\}$ and $\{Y^\prime_t: t\ge 0\}$
are two continuous solutions of (\ref{3.8}). Fix $R \ge 1$ and let
$\tau = \inf\{t\ge0 : \<\phi_p, Y_t\> \ge R$ or $\<\phi_p,
Y^\prime_t\> \ge R\}$. Observe that $q(Y_t,x) 1_{\{t\le\tau\}}$ and
$q(Y^\prime_t,x) 1_{\{t\le\tau\}}$ are predictable and define
 \beqnn
X_t^* = \int_0^{t\wedge\tau}\int_{\mbb{R}}\int_0^{q(Y_s,a) \wedge
q(Y^\prime_s,a)} \int_{W_0} w(t-s)\delta_{x_{s,a}(t)}
N_1(ds,da,du,dw)
 \eeqnn
and
 \beqnn
Y_t^* = \int_0^{t\wedge\tau}\int_{\mbb{R}}\int_0^{q(Y_s,a) \vee
q(Y^\prime_s,a)} \int_{W_0} w(t-s)\delta_{x_{s,a}(t)}
N_1(ds,da,du,dw).
 \eeqnn
By Proposition~\ref{p2.2} we have
 \beqlb\label{3.9}
\mbf{E}[\<\phi,X_t^*\>] = \int_0^t \mbf{E}[\<(q(Y_s,\cdot) \land
q(Y^\prime_s,\cdot)) P_{t-s}\phi,m\> 1_{\{s\le\tau\}}]ds
 \eeqlb
and
 \beqlb\label{3.10}
\mbf{E}[\<\phi,Y_t^*\>] = \int_0^t \mbf{E}[\<(q(Y_s,\cdot) \vee
q(Y^\prime_s,\cdot)) P_{t-s}\phi,m\> 1_{\{s\le\tau\}}]ds.
 \eeqlb
Let $c\ge0$ be a constant such that $|\phi_p^{\prime\prime}| \le
2c\phi_p$. Then we have
 \beqnn
\frac{d}{dt}P_t\phi_p(x) = \frac{1}{2}P_t\phi_p^{\prime\prime}(x)
\le cP_t\phi_p(x).
 \eeqnn
By a comparison theorem, we get $P_t\phi_p \le e^{ct}\phi_p$ for
every $t\ge0$. Observe also that $\|Y^\prime_{t\wedge\tau} -
Y_{t\wedge\tau}\|_p \le \<\phi_p, Y^*_t\> - \<\phi_p, X^*_t\>$. Then
(\ref{3.9}) and (\ref{3.10}) imply that
 \beqlb\label{3.11}
\mbf{E}[\|Y^\prime_{t\wedge \tau} - Y_{t\wedge\tau}\|_p]
 \ar=\ar
\int_0^t \mbf{E}[\<|q(Y^\prime_s,\cdot) - q(Y_s,\cdot)|
P_{t-s}\phi_p,m\>1_{\{s\le\tau\}}]ds  \nnm  \\
 \ar\le\ar
e^{ct}\int_0^t \mbf{E}[\<|q(Y^\prime_s,\cdot) -
q(Y_s,\cdot)|\phi_p,m\>1_{\{s\le\tau\}}]ds  \nnm \\
 \ar\le\ar
L_Re^{ct}\int_0^t \mbf{E}[\|Y^\prime_{s\wedge\tau} -
Y_{s\wedge\tau}\|_p]ds.
 \eeqlb
By Gronwall's inequality we conclude $\mbf{E}[\|Y^\prime_{t\wedge
\tau} - Y_{t\wedge\tau}\|_p] = 0$. Since $R\ge0$ can be arbitrary,
that gives the uniqueness of solution. To show the existence of a
solution, we first assume (\ref{3.3}) holds with $L_R$ replaced by a
universal constant $L$ independent of $R$. We define the sequence
continuous processes $\{Y_t^{(n)}: t\ge0\}$ by setting $Y_t^{(0)} =
\sum_{i=1}^\infty \xi_i(t) \delta_{x_{0,a_i}(t)}$ and
 \beqlb\label{3.12}
Y_t^{(n)} = Y_t^{(0)} + \int_0^t\int_{\mbb{R}}
\int_0^{q(Y_s^{(n-1)},a)} \int_{W_0} w(t-s)
\delta_{x_{s,a}(t)}N_1(ds,da,du,dw)
 \eeqlb
for $n\ge1$. By a reasoning as in (\ref{3.11}) we see that
 \beqnn
\mbf{E}[\|Y^{(n)}_t - Y^{(n-1)}_t\|_p]
 \le
Le^{ct}\int_0^t \mbf{E}[\|Y^{(n-1)}_s - Y^{(n-2)}_s\|_p]ds.
 \eeqnn
By a standard argument, one sees there is a predictable process
$\{Y_t: t\ge0\}$ such that $\lim_{n\to \infty} \mbf{E}[\|Y_t -
Y^{(n)}_t\|_p] = 0$ uniformly on each bounded interval of $t\ge0$.
By Theorem~\ref{t2.1} we see that
 \beqlb\label{3.13}
\tilde Y_t = Y_t^{(0)} + \int_0^t\int_{\mbb{R}} \int_0^{q(Y_s,a)}
\int_{W_0} w(t-s) \delta_{x_{s,a}(t)}N_1(ds,da,du,dw)
 \eeqlb
defines a continuous process $\{\tilde Y_t: t\ge0\}$ on
$M_p(\mbb{R})$. Based on (\ref{3.12}) and (\ref{3.13}) we may follow
the calculations in (\ref{3.11}) to obtain
 \beqnn
\mbf{E}[\|Y^{(n)}_t - \tilde Y_t\|_p]
 \le
Le^{ct}\int_0^t \mbf{E}[\|Y^{(n-1)}_s - Y_s\|_p]ds.
 \eeqnn
Letting $n\to \infty$ we see that $Y_t = \tilde Y_t$ a.s.\ for every
$t\ge0$. That proves the existence of a solution of (\ref{3.8}) for
a universal constant $L$. The extension of the existence to the
general condition (\ref{3.3}) is a standard localization argument.
The strong Markov property of $\{Y_t: t\ge 0\}$ follows from the
uniqueness of the solution. \qed

To consider a general initial value $\mu \in M_p^a(\mbb{R})$ let
$\{W(ds,dy)\}$, $\{N_0(da,dw)\}$ and $\{N_1(ds,da$, $du,dw)\}$ be
given as in (2.a), (2.c) and (2.d). Let $(\mcr{F}_t)_{t\ge 0}$ be
the filtration generated by (\ref{2.14}) and (\ref{2.15}). Let
$\{X_t: t\ge 0\}$ be given by (\ref{2.16}). By a \textit{solution}
of the stochastic equation
 \beqlb\label{3.14}
Y_t = X_t + \int_0^t\int_{\mbb{R}}\int_0^{q(Y_s,a)}\int_{W_0}
w(t-s) \delta_{x_{s,a}(t)}N_1(ds,da,du,dw), \qquad t\ge0,
 \eeqlb
we mean an $(\mcr{F}_t)$-adapted continuous process $\{Y_t:
t\ge0\}$ in $M_p(\mbb{R})$ that satisfies (\ref{3.14}). By
arguments similar to those in the proof of Theorem~\ref{t3.1} we
have

\btheorem\label{t3.2} There is a unique solution $\{Y_t: t\ge0\}$
of (\ref{3.14}), which also solves the martingale problem given by
(\ref{3.4}) and (\ref{3.5}). \etheorem


\section{Existence and continuity of local times}

\setcounter{equation}{0}

In this section, we prove that the occupation time of an immigration
superprocess $\{Y_t: t\ge 0\}$ has a H\"older continuous density
field $\{z(b,t): b\in \mbb{R}, t\ge 0\}$. As in the classical case
the two parameter process $\{z(b,t): b\in \mbb{R}, t\ge 0\}$ can be
interpreted as the local time of $\{Y_t: t\ge 0\}$. We shall need
the following two lemmas.

\blemma\label{l4.1} Let $\{B_j(\cdot): j=1,\cdots,n\}$ be a family
of Brownian motions, $\{\xi_j(\cdot): j=1,\cdots,n\}$ be a family
of independent $\sigma$-branching diffusions and
$\{\alpha_j(\cdot): j=1,\cdots,n\}$ be a family of bounded
processes. Suppose that the two families $\{B_j(\cdot):
j=1,\cdots,n\}$ and $\{\xi_j(\cdot): j=1,\cdots,n\}$ are
independent of each other. Let $l_j(\cdot,\cdot)$ denote the local
time of $B_j(\cdot)$. Then to each integer $k\ge1$ there
corresponds a constant $C_k\ge0$ such that
 \beqlb\label{4.1}
 \ar \ar\mbf{E}\bigg\{\bigg[\sum_{j=1}^n\int_{r_j}^{t_j} \xi_j(s)
\alpha_j(s) l_j(b,ds)\bigg]^{2k}\bigg\}  \nnm \\
 \ar \ar\qquad
\le C_k(t-r)^k\bigg\{\mbf{E}\bigg[\bigg(\sum_{j=1}^n
\xi_j(t_j)\bigg)^{2k}\bigg] + \mbf{E}\bigg[\bigg(\sum_{j=1}^n
\int_{r_j}^{t_j}\xi_j(s)ds\bigg)^k\bigg]\bigg\}
 \eeqlb
for any intervals $[r_j,t_j] \subset [r,t]$, $j=1,\cdots,n$.
\elemma

\proof In this and the following proofs, $C_k$ will denote
positive constants that may change values from line to line. By
integration by parts formula, we have
 \beqlb\label{4.2}
\int_{r_j}^{t_j} \xi_j(s)\alpha_j(s) l_j(b,ds)
 \ar=\ar
\xi_j(t_j) \int_{r_j}^{t_j} \alpha_j(s) l_j(b,ds) \nnm \\
 \ar \ar\qquad
- \int_{r_j}^{t_j} d\xi_j(s)\int_{r_j}^s \alpha_j(s) l_j(b,ds).
 \eeqlb
By H\"older's inequality and the boundedness of $\{\alpha_j
(\cdot): j=1,\cdots,n\}$,
 \beqlb\label{4.3}
 \ar\ar\mbf{E}\bigg\{\bigg[\sum_{j=1}^n\xi_j(t_j) \int_{r_j}^{t_j}
\alpha_j(s) l_j(b,ds)\bigg]^{2k}\bigg\} \nnm \\
 \ar\ar\qquad
\le \mbf{E}\bigg\{\bigg[\sum_{j=1}^n\xi_j(t_j)\bigg]^{2k-1}
\sum_{j=1}^n \xi_j(t_j) \bigg[\int_{r_j}^{t_j}
\alpha_j(s) l_j(b,ds)\bigg]^{2k}\bigg\} \nnm \\
 \ar\ar\qquad
\le C_k\mbf{E}\bigg\{\bigg[\sum_{j=1}^n\xi_j(t_j)\bigg]^{2k-1}
\sum_{j=1}^n \xi_j(t_j)l_j(b,[r_j,t_j])^{2k}\bigg\} \nnm \\
 \ar\ar\qquad
\le C_k(t-r)^k\mbf{E}\bigg\{\bigg[\sum_{j=1}^n
\xi_j(t_j)\bigg]^{2k}\bigg\},
 \eeqlb
where we also used the estimate $\mbf{E}[l_j(b, [r_j,t_j])^{2k}]
\le C_k(t_j-r_j)^k$ and the independence of $\xi_j(\cdot)$ and
$l_j(\cdot,\cdot)$ for the last inequality. Similarly, by Doob's
martingale inequality we get
 \beqlb\label{4.4}
 \ar\ar \mbf{E}\bigg\{\bigg[\sum_{j=1}^n \int_{r_j}^{t_j} d\xi_j(s)
\int_{r_j}^s \alpha_j(s) l_j(b,ds)\bigg]^{2k}\bigg\}  \nnm \\
 \ar\ar\qquad
\le C_k\mbf{E}\bigg\{\bigg[\sum_{j=1}^n \int_{r_j}^{t_j}
\bigg(\int_{r_j}^s \alpha_j(s) l_j(b,ds)\bigg)^2\xi_j(s)ds
\bigg]^k\bigg\}  \nnm \\
 \ar\ar\qquad
\le C_k\mbf{E}\bigg\{\bigg[\sum_{j=1}^n \int_{r_j}^{t_j}
l_j(b,[r_j,s])^2\xi_j(s)ds\bigg]^k\bigg\}  \nnm \\
 \ar\ar\qquad
\le C_k\mbf{E}\bigg\{\bigg[\sum_{j=1}^n \int_{r_j}^{t_j}
\xi_j(s)ds \bigg]^{k-1} \sum_{j=1}^n\int_{r_j}^{t_j}
l_j(b,[r_j,s])^{2k}\xi_j(s)ds\bigg\} \nnm \\
 \ar\ar\qquad
\le C_k(t-r)^k\mbf{E}\bigg\{\bigg[\sum_{j=1}^n \int_{r_j}^{t_j}
\xi_j(s) ds\bigg]^k\bigg\}.
 \eeqlb
Then the desired estimate follows from (\ref{4.2}), (\ref{4.3})
and (\ref{4.4}). \qed

\blemma\label{l4.2} Under the assumptions of Lemma~\ref{l4.1}, to
each integer $k\ge1$ there corresponds a constant $C_k\ge0$ such
that
 \beqlb\label{4.5}
 \ar\ar \mbf{E}\bigg\{\bigg|\sum_{j=1}^n\bigg[\int_{r_j}^{t_j} \xi_j(s)
\alpha_j(s)l_j(b_1,ds) - \int_{r_j}^{t_j} \xi_j(s)\alpha_j(s)
l_j(b_2,ds)\bigg]\bigg|^{2k}\bigg\}  \nnm \\
 \ar\ar\qquad
\le C_k(b_1-b_2)^k\bigg\{\mbf{E}\bigg[\bigg(\sum_{j=1}^n
\xi_j(t_j)\bigg)^{2k}\bigg] + \mbf{E}\bigg[\bigg(\sum_{j=1}^n
\int_{r_j}^{t_j}\xi_j(s)ds\bigg)^k\bigg]\bigg\}
 \eeqlb
for any $b_1,b_2\in\mbb{R}$ and any intervals $[r_j,t_j] \subset
[r,t]$, $j=1,\cdots,n$. \elemma

\proof By a formula at page 211 of Ref.~\cite{RY91} it is easy to
show that
 \beqnn
\mbf{E}\big\{[l_j(b_1, [r_j,t_j]) - l_j (b_2,
[r_j,t_j])\big]^{2k}\big\}
 \le
C_k(b_1-b_2)^k.
 \eeqnn
We can apply (\ref{4.2}) to the two integrals on the left hand side
of (\ref{4.5}). Then the result is obtained by similar estimates as
in the proof of Lemma~\ref{l4.1}. \qed

\btheorem\label{t4.1} Let $\eta(\cdot,\cdot,\cdot)\in
L^2_+(\mcr{P})$ and suppose there is a deterministic increasing
function $\bar \eta(\cdot)$ on $[0,\infty)$ such that
$\eta(t,a,\omega)\le \bar \eta(t)$ for all $(t,a,\omega)$. Let
$\{Z_t: t\ge0\}$ denote the occupation time of the process $\{Y_t:
t\ge0\}$ defined by (\ref{2.17}). Then $Z_t$ is a.s.\ absolutely
continuous and the corresponding local time is given by
 \beqlb\label{4.6}
z(b,t)
 \ar=\ar
  \int_{\mbb{R}}\int_{W_0}N_0(da,dw)\int_0^t w(u)
l_{0,a}(b,du)  \nnm  \\
 \ar \ar
+
\int_0^t\int_{\mbb{R}}\int_0^{\eta(s,a)}\int_{W_0}N_1(ds,da,du,dw)
\int_s^t w(v-s)l_{s,a}(b,dv),
 \eeqlb
where $\{l_{s,a}(b,u): u\ge s\}$ is the local time of
$\{x_{s,a}(u): u\ge s\}$ at $b\in \mbb{R}$. Moreover, the two
parameter process $\{z(b,t): b\in \mbb{R}, t\ge 0\}$ has a version
which is H\"older continuous of order $\alpha$ for every $\alpha <
1/2$. \etheorem

\proof The existence of the local time follows by (\ref{2.17}) and
Fubini's theorem. Let $z_0(b,t)$ and $z_1(b,t)$ denote respectively
the first and the second terms on the right hand side of
(\ref{4.6}). We shall only give the proof of the continuity result
for $z_1(\cdot,\cdot)$ since the proof for $z_0(\cdot,\cdot)$ is
similar.

\textit{Step~1.} Let us consider the special case with $\mu$ and
$m\in M(\mbb{R})$. Observe that $z_1(b,0) = 0$ for every $b\in
\mbb{R}$. For any $t>r>0$ and $b\in \mbb{R}$ we have
 \beqnn
z_1(b,t)-z_1(b,r)
 \ar=\ar \int_r^t\int_{\mbb{R}}\int_0^{\eta(s,a)}\int_{W_0}N_1(ds,da,du,dw)
\int_s^t w(v-s)l_{s,a}(b,dv) \nnm \\
 \ar \ar
+
\int_0^r\int_{\mbb{R}}\int_0^{\eta(s,a)}\int_{W_0}N_1(ds,da,du,dw)
\int_r^t w(v-s)l_{s,a}(b,dv).
 \eeqnn
For $\epsilon>0$ let $W_\epsilon = \{w\in W_0: w(\epsilon)>0\}$.
Then $N_1([0,t]\times \mbb{R} \times [0,\bar \eta(t)] \times
W_\epsilon)$ is a.s.\ finite for every $t\ge0$. Let $S(\bar q,N_1)$
denote the intersection of the support of $N_1(ds,da,du,dw)$ with
the set $\{(s,a,u,w): s\ge0, a\in \mbb{R}, 0\le u\le \bar \eta(s),
w\in W_\epsilon\}$. We can enumerate $S(\bar q,N_1)$ into a sequence
$\{(s_j,a_j,u_j,w_j): j=1,2,\cdots\}$ so that $0< s_1< s_2< \cdots$.
As in the proof of Lemma~3.2 in Ref.~\cite{DL03} one sees that,
given $\{(s_j,a_j,u_j,w_j(\epsilon)): j=1,2,\cdots\}$ each
$\{w_j(u): u\ge \epsilon\}$ is a $\sigma$-branching diffusion
independent of the white noise $\{W(ds,dy)\}$. For any integer
$k\ge1$ we may use Lemma~\ref{l4.1} and Fatou's lemma to see that
 \beqnn
 \ar\ar\mbf{E}\bigg\{\bigg[\int_r^t\int_{\mbb{R}}\int_0^{\eta(s,a)}\int_{W_0}
N_1(ds,da,du,dw)\int_s^t w(v-s)l_{s,a}(b,dv)\bigg]^{2k}\bigg\}
 \nnm  \\
 \ar\ar\qquad
\le \liminf_{\epsilon\to0}\mbf{E}\bigg\{\bigg[\sum_{r<s_j\le t}
\int_{s_j+\epsilon}^t w_j(v-s_j)
l_{s_j,a_j}(b,dv)\bigg]^{2k}\bigg\} \nnm  \\
 \ar\ar\qquad
\le C_k(t-r)^k\liminf_{\epsilon\to0} \bigg\{\mbf{E}
\bigg[\bigg(\sum_{r<s_j\le t}w_j(t-s_j)\bigg)^{2k}\bigg]
 \nnm  \\
 \ar\ar\qquad\qquad
+ \mbf{E}\bigg[\bigg(\sum_{r<s_j\le t}\int_{\epsilon+s_j}^t
w_j(v-s_j)dv\bigg)^k\bigg]\bigg\} \nnm  \\
 \ar\ar\qquad
\le C_k(t-r)^k\bigg\{\mbf{E} \bigg[\bigg(\int_r^t\int_{\mbb{R}}
\int_0^{\bar \eta(t)}\int_{W_0}w(t-s)N_1(ds,da,du,dw)\bigg)^{2k}
\bigg] \nnm  \\
 \ar\ar\qquad\qquad
+ \mbf{E}\bigg[\bigg(\int_r^t\int_{\mbb{R}}\int_{W_0} \int_0^{\bar
\eta(t)} N_1(ds,da,du,dw)\int_s^t w(v-s)dv\bigg)^k
\bigg]\bigg\} \nnm  \\
 \ar\ar\qquad
\le C_k(t-r)^k\Big\{\mbf{E}\big[\<1,\bar Y_t\>^{2k}\big] +
\mbf{E}\big[\<1,\bar Z_t\>^k\big]\Big\},
 \eeqnn
where $\{\bar Y_t\}$ is the process defined by (\ref{2.17}) from
$\{\bar \eta(t)\}$, and $\{\bar Z_t\}$ is the corresponding
occupation time process. By a similar reasoning as the above we have
 \beqnn
\ar\ar
\mbf{E}\bigg\{\bigg[\int_0^r\int_{\mbb{R}}\int_0^{\eta(s,a)}
\int_{W_0} N_1(ds,da,du,dw)\int_r^t w(u-s)l_{s,a}(b,du)\bigg]^{2k}
\bigg\} \nnm  \\
 \ar\ar\qquad
\le \liminf_{\epsilon\to 0}\mbf{E}\bigg\{\bigg[\sum_{0<s_j\le r}
\int_{r\vee (\epsilon+s_j)}^t w_j(u-s_j)
l_{s_j,a_j}(b,du)\bigg]^{2k}\bigg\} \nnm  \\
 \ar\ar\qquad
\le C_k(t-r)^k\liminf_{\epsilon\to 0} \bigg\{\mbf{E}
\bigg[\bigg(\sum_{0<s_j\le r}w_j(t-s_j)\bigg)^{2k}\bigg] \nnm  \\
 \ar\ar\qquad\qquad
+ \mbf{E}\bigg[\bigg(\sum_{0<s_j\le r} \int_{r\vee
(\epsilon+s_j)}^t w_j(u-s_j)du\bigg)^k\bigg]\bigg\} \nnm  \\
 \ar\ar\qquad
\le C_k(t-r)^k\bigg\{\mbf{E} \bigg[\bigg(\int_0^r\int_{\mbb{R}}
\int_0^{\bar \eta(t)}\int_{W_0}w(t-s)N_1(ds,da,du,dw)
\bigg)^{2k}\bigg] \nnm  \\
 \ar\ar\qquad\qquad
+ \mbf{E}\bigg[\bigg(\int_0^r\int_{\mbb{R}}\int_0^{\bar \eta(t)}
\int_{W_0}N_1(ds,da,du,dw)\int_r^t w(v-s)dv\bigg)^k
\bigg]\bigg\} \nnm \\
 \ar\ar\qquad
\le C_k(t-r)^k\Big\{\mbf{E}\big[\<1,\bar Y_t\>^{2k}\big] +
\mbf{E}\big[\<1,\bar Z_t\>^k\big]\Big\}.
 \eeqnn
Thus we have
 \beqlb\label{4.7}
\mbf{E}[|z_1(b,t)-z_1(b,r)|^{2k}] \le C_k(t-r)^k
\Big\{\mbf{E}\big[\<1, \bar Y_t\>^{2k}\big] + \mbf{E}\big[\<1,\bar
Z_t\>^k\big]\Big\}.
 \eeqlb
For $t>0$ and $b_1, b_2\in \mbb{R}$ we can use Lemma~\ref{l4.2}
and similar arguments as the above to show
 \beqlb\label{4.8}
\mbf{E}[|z_1(b_1,t)-z_1(b_2,t)|^{2k}] \le C_k(b_1-b_2)^k
\Big\{\mbf{E}\big[\<1,\bar Y_t\>^{2k}\big] + \mbf{E}\big[\<1,\bar
Z_t\>^k\big]\Big\}.
 \eeqlb
Since $\mu$ and $m\in M(\mbb{R})$ are finite measures and $\bar
\eta(t)$ is a locally bounded function of $t\ge0$, it is not hard to
show that $\mbf{E}[\<1, \bar Y_s\>^{2k} + \<1, \bar Z_s\>^k]$ is
locally bounded in $t\ge0$. Then (\ref{4.7}) and (\ref{4.8}) imply
that $z_1(\cdot,\cdot)$ has a H\"older continuous version of order
$\alpha$ for every $\alpha < 1/2$; see e.g.\ page 273 of
Ref.~\cite{W86}.

\textit{Step~2.} Now we consider the general case with $\mu$ and
$m\in M_p(\mbb{R})$. Let $L$ and $T$ be fixed positive constants.
For any integer $n\ge 1$ let $a_n>0$ be sufficiently large so that
 \beqlb\label{4.9}
\mbf{P}\Big\{\sup_{0\le s\le T} |x_{0,0}(s)| \ge a_n\Big\} \le
\frac{1}{n}.
 \eeqlb
Let $b_n=L+a_n$ and $c_n=L+2a_n$. Let $A_n= \{$there exists $0\le
s\le T$ such that $x_{0,b_n}(s) = L$ or $=c_n\}$ and $B_n= \{$there
exists $0\le s\le T$ such that $x_{0,-b_n}(s) = -L$ or $=-c_n\}$. As
observed in Ref.~\cite{W98}, any two solutions of (\ref{2.12})
started from different locations never collide. Then on the event
$(A_n\cup B_n)^c$, for any $|a|\ge c_n$ and $0\le s\le u\le T$ we
have $x_{s,a}(u)> x_{0,b_n}(u)> L$ or $x_{s,a}(u)< x_{0,-b_n}(u)<
-L$, and so $l_{s,a}(b,u)=0$ whenever $|b|\le L$. It follows that
 \beqnn
z(b,t)
 \ar=\ar
\int_{[-c_n,c_n]}\int_{W_0}N_0(da,dw)\int_0^t w(u)
l_{0,a}(b,du)  \nnm  \\
 \ar \ar
+ \int_0^t\int_{[-c_n,c_n]}\int_0^{\eta(s,a)}\int_{W_0}
N_1(ds,da,du,dw) \int_s^t w(v-s)l_{s,a}(b,dv)
 \eeqnn
for $(b,t) \in [-L,L]\times [0,T]$ on $(A_n\cup B_n)^c$. By Step~1,
for any $\alpha < 1/2$ we have a modification of $z(b,t)$ on
$(A_n\cup B_n)^c$ that is H\"older continuous of order $\alpha$ in
$(b,t) \in (-L,L)\times [0,T)$. In view of (\ref{4.9}) we have
$\mbf{P}(A_n) = \mbf{P}(B_n) =1/n$. Then we can modify $z(b,t)$ on
the whole space $\Omega$ so the process becomes H\"older continuous
of order $\alpha$ in $(b,t) \in (-L,L)\times [0,T)$. Since $L\ge0$
and $T\ge0$ can be arbitrary, we have the desired result. \qed

\btheorem\label{t4.2} Let $q(\cdot,\cdot)$ be a Borel function on
$M_p(\mbb{R}) \times \mbb{R}$ satisfying the conditions in the last
section. In stead of (\ref{3.2}), we assume the stronger condition
 \beqlb\label{4.9a}
q(\mu,a)\le K(1+\<\phi_p,\mu\>^2), \qquad  a \in \mbb{R}, \mu\in
M_p(\mbb{R}).
 \eeqlb
Let $\{Z_t: t\ge0\}$ denote the occupation time of the ISDSM defined
by (\ref{3.14}). Then $Z_t$ is a.s.\ absolutely continuous and the
corresponding local time is given by
 \beqlb\label{4.10}
z(b,t)
 \ar=\ar
  \int_{\mbb{R}}\int_{W_0}N_0(da,dw)\int_0^t w(u)
l_{0,a}(b,du)  \nnm  \\
 \ar \ar
+ \int_0^t\int_{\mbb{R}}\int_0^{q(Y_s,a)}\int_{W_0}
N_1(ds,da,du,dw) \int_s^t w(v-s)l_{s,a}(b,dv).
 \eeqlb
Moreover, the two parameter process $\{z(b,t): b\in \mbb{R}, t\ge
0\}$ has a version which is H\"older continuous of order $\alpha$
for every $\alpha < 1/2$. \etheorem

\proof Let $\{Y_t^{(n)}: t\ge0\}$ be defined by the right hand
side of (\ref{3.13}) with $q(Y_s,a)$ replaced by $\eta_n(s,a) :=
n\land q(Y_s,a)$. {From} Theorem~\ref{t4.1} we know that the
occupation time of $\{Y_t^{(n)}: t\ge 0\}$ is absolutely
continuous with density given by
 \beqnn
z_n(b,t)
 \ar=\ar
  \int_{\mbb{R}}\int_{W_0}N_0(da,dw)\int_0^t w(u)
l_{0,a}(b,du)  \nnm  \\
 \ar \ar
+ \int_0^t\int_{\mbb{R}}\int_0^{\eta_n(s,a)}\int_{W_0}
N_1(ds,da,du,dw) \int_s^t w(v-s)l_{s,a}(b,dv)
 \eeqnn
which has a H\"older continuous version of order $\alpha$. The
same conclusion is clearly true for the occupation time of $\{Y_t:
t\ge 0\}$ on the event $\{\sup_{0\le s\le t}q(Y_s,a)\le n\}$. By
Chebyshev's inequality and the assumption on $q(\cdot,\cdot)$ is
easy to show that
 \beqnn
\mbf{P}\Big\{\sup_{0\le s\le t}q(Y_s,a)\ge n\Big\}
 \le
\frac{K}{n}\mbf{E}\Big[1+\sup_{0\le s\le t}\<\phi_p,Y_s\>^2\Big].
 \eeqnn
By Proposition~\ref{p3.1}, the right hand side tends to zero as
$n\to \infty$. Then we have the desired result. \qed


\section{Scaling limit theorems}

\setcounter{equation}{0}

Scaling limit theorems of SDSM without immigration were investigated
in Refs.~\cite{DLZ04,LWX04b}. A direct generalization of the limit
theorem of Ref.~\cite{DLZ04} to the ISDSM was given in
Ref.~\cite{D06}. In this section, we prove two scaling limit
theorems for the ISDSM. The limit processes obtained here are
different from those in the previous work. As consequences, we also
obtain scaling limits for the corresponding local time. For
simplicity we assume $Y_0=0$ and focus on the influence of the
immigration.

Let $\bar\mbb{R} = \mbb{R} \cup\{\infty\}$ be the one-point
compactification of the real line and let $M(\bar\mbb{R})$ be the
space of Borel measures on $\bar\mbb{R}$. We fix a metric on
$M(\bar\mbb{R})$ compatible with the weak convergence and regard
$M(\mbb{R})$ as a subspace of $M(\bar\mbb{R})$ comprising measures
supported by $\mbb{R}$. Then a metric on $M_p(\mbb{R})$ can be
defined through the isomorphism
 \beqlb\label{5.0}
\Phi_p: \mu(dx) \mapsto \phi_p(x)\mu(dx)
 \eeqlb
between $M_p(\mbb{R})$ and $M(\mbb{R})$. Let $C([0,\infty),
M_p(\bar\mbb{R}))$ be the space of continuous paths from
$[0,\infty)$ to $M_p(\bar\mbb{R})$ endowed with the topology of
locally uniform convergence.

We first assume $p>1$ so the Lebesgue measure $\lambda$ belongs to
$M_p(\mbb{R})$. Let $q(\cdot,\cdot)$ be a bounded Borel function on
$M_p(\mbb{R})\times \mbb{R}$ satisfying the local Lipschitz
condition (\ref{3.3}). Suppose that $\{W(ds,dy)\}$ and
$\{N_1(ds,da$, $du,dw)\}$ are given as in (2.a) and (2.d) with $m$
replaced by $\lambda$. Let $\{x_{r,a}(t): t\ge r\}$ be defined by
(\ref{2.12}) and let $\{Y_t: t\ge0\}$ be the solution of
 \beqlb\label{5.1}
Y_t = \int_0^t\int_{\mbb{R}}\int_0^{q(Y_s,a)}\int_{W_0} w(t-s)
\delta_{x_{s,a}(t)}N_1(ds,da,du,dw).
 \eeqlb
For any integer $k\ge1$ let $Y_t^k(dx) = k^{-2}Y_{k^2t}(kdx)$.

\btheorem\label{t5.1} Suppose that $q(\nu,a) \to q(\infty)$ as
$|a| \to\infty$ for all $\mu \in M_p(\mbb{R})$. Then, as $k\to
\infty$, $\{k^{-1}Y_t^k: t\ge0\}$ converges to $\{q(\infty)
t\lambda: t\ge0\}$ in probability on $C([0,\infty),
M_p(\mbb{R}))$. \etheorem

\proof In this and the following proofs, we write ``$=_d$'' for
the equivalence in distribution of two processes. We shall also
consider some new Poisson random measures which might be defined
on some extensions of the original probability space.

\textit{Step~1.} Let $\phi\in C^2_p(\mbb{R})$. {From} (\ref{5.1}) we
have
 \beqlb\label{5.2}
\<\phi,Y_t^k\>
 =
\int_0^t\int_{\mbb{R}}\int_0^{q(Y_{k^2s},kb)}\int_{W_0} w_k(t-s)
\phi(x^k_{s,b}(t)) N_1(k^2ds,k db,du,dw),
 \eeqlb
where $w_k(t-s) = k^{-2}w(k^2(t-s))$ and $\{x_{s,b}^k(t): t\ge
s\}$ is the unique solution of
 \beqnn
x(t) = b + \int_s^t\int_{\mbb{R}} h_k(y-x(u)) k^{-3/2} W(k^2du,k
dy), \qquad t\ge s
 \eeqnn
with $h_k(z) = \sqrt{k}h(k z)$. It is easy to see that $k^{-3/2}
W(k^2du,k dy)$ is a white noise based on the Lebesgue measure and
$N_1(k^2ds,k db,du,dw)$ is Poisson random measure with intensity
$k^3dsdbdu \mbf{Q}_\kappa(dw)$. By the scaling property of the
$\sigma$-branching diffusion it is easy to check that $\{w_k(t):
t\ge0\}$ under $k^2 \mbf{Q}_\kappa$ has the same law as $\{w(t):
t\ge0\}$ under $\mbf{Q}_\kappa$. {From} (\ref{5.2}) we get
 \beqlb\label{5.3}
\<\phi,Y_t^k\>
 =_d
\int_0^t\int_{\mbb{R}}\int_0^{q(Y_{k^2s},kb)}\int_{W_0} w(t-s)
\phi(x^k_{s,b}(t)) N_k(ds,db,du,dw),
 \eeqlb
where $N_k(ds,db,du,dw)$ is Poisson random measure with intensity
$kdsdbdu\mbf{Q}_\kappa(dw)$. Let $\psi_k = k^{-1}\phi$. Since
$\mbf{Q}_\kappa[w(t)]=1$ for every $t>0$, from (\ref{5.3}) we have
 \beqnn
\mbf{E}[\<\psi_k,Y_t^k\>]
 \ar=\ar
\mbf{E}\bigg[\int_0^tds\int_{\mbb{R}}db \int_{W_0} q(Y_{k^2s},kb)
w(t-s) \phi(x^k_{s,b}(t))\mbf{Q}_\kappa(dw)\bigg]  \\
 \ar=\ar
\int_0^tds\int_{\mbb{R}}\mbf{E}[q(Y_{k^2s},kb)\phi(x^k_{s,b}(t))]db  \\
 \ar=\ar
\int_0^tds\int_{\mbb{R}}db\int_{\mbb{R}}\mbf{E}[q(Y_{k^2s},kb)|
x^k_{s,b}(t)=z]\phi(z)g_{\rho(0)(t-s)}(b,z)dz,
 \eeqnn
where $g_u(b,z)$ denotes the density of the heat kernel. By
dominated convergence we see that
 \beqnn
\lim_{k\to \infty} \mbf{E}[\<\psi_k, Y_t^k\>]
 =
\int_0^tds\int_{\mbb{R}}db\int_{\mbb{R}}
q(\infty)\phi(z)g_{\rho(0)(t-s)}(b,z)dz
 =
t q(\infty)\<\phi,\lambda\>.
 \eeqnn
For every $t>0$ we have $\mbf{Q}_\kappa[w(t)^2] = \sigma t$; see
e.g.\ Lemma~3.1 of Ref.~\cite{S90}. It follows that
 \beqnn
\ar\ar\mbf{E}\Big[\Big|\<\psi_k,Y_t^k\> -
\mbf{E}[\<\psi_k,Y_t^k\>]\Big|^2\Big]  \\
 \ar\ar\qquad
= k^{-1}\mbf{E}\bigg[\int_0^tds\int_{\mbb{R}}db \int_{W_0}
q(Y_{k^2s},kb)
w(t-s)^2\phi(x^k_{s,b}(t))^2\mbf{Q}_\kappa(dw)\bigg]  \\
 \ar\ar\qquad
\le  k^{-1}\sigma \|q\|\int_0^t(t-s)ds\int_{\mbb{R}}\mbf{E}
\big[\phi(x^k_{s,b}(t))^2\big]db  \\
 \ar\ar\qquad
\le k^{-1}\sigma \|q\| \int_0^t(t-s)ds\int_{\mbb{R}}\phi(b)^2db,
 \eeqnn
which tends to zero as $k\to \infty$. That proves $\lim_{k\to
\infty} \<\psi_k, Y_t^k\> = q(\infty)t\<\phi, \lambda\>$ in
probability.

\textit{Step~2.} {From} (\ref{5.3}) and Theorem~\ref{t2.1} we see
that
 \beqlb\label{5.4}
M_t^k(\phi) := \<\phi, Y_t^k\> - \frac{1}{2} \rho(0)\int_0^t
\<\phi^{\prime\prime}, Y_s^k\> ds - k\int_0^t \<q(Y_{k^2s},
k\cdot)\phi,\lambda\>ds
 \eeqlb
is a continuous martingale with quadratic variation process
 \beqlb\label{5.5}
\<M^k(\phi)\>_t = \int_0^t\<\sigma\phi^2, Y_s^k\>ds + \int_0^t
ds\int_{\mbb{R}^2} \rho(k(x-y)) \phi^\prime(x)\phi^\prime(y)
Y_s^k(dx)Y_s^k(dy).
 \eeqlb
Let $\{\tau_k\}$ be a bounded family of stopping times. By
(\ref{5.4}) and (\ref{5.5}) we have
 \beqnn
\mbf{E}\big[\big|\<\phi, Y_{\tau_k+t}^k\> - \<\phi,
Y_{\tau_k}^k\>\big|^2\big]
 \ar\le\ar
3\mbf{E}\big[|M_{\tau_k+t}^k(\phi) - M_{\tau_k}^k(\phi)|^2] +
\frac{3}{2}\rho(0)t\int_0^t \mbf{E}\big[\<\phi^{\prime\prime},
Y_{\tau_k+s}^k\>^2\big] ds
 \\
 \ar \ar
+\, 3k^2t\int_0^t \mbf{E}\big[\<q(Y_{k^2(\tau_k+s)},
k\cdot)\phi,\lambda\>^2\big]ds,
 \eeqnn
where
 \beqnn
\mbf{E}\big[|M_{\tau_k+t}^k(\phi) - M_{\tau_k}^k(\phi)|^2]
 \ar\le\ar
\sigma\int_0^t \mbf{E}\big[\<\phi^2, Y_{\tau_k+s}^k\>\big] ds +
\|\rho\|\int_0^t \mbf{E}\big[\<\phi^\prime, Y_{\tau_k +
s}^k\>^2\big]ds,
 \eeqnn
By Proposition~\ref{p2.1} it is easy to show that
 \beqnn
\lim_{t\to0}\sup_{k\ge1}\mbf{E}\big[\big|\<\psi_k, Y_{\tau_k+t}^k\>
- \<\psi_k, Y_{\tau_k}^k\>\big|^2\big] = 0.
 \eeqnn
Another application of Proposition~\ref{p2.1} shows
 \beqnn
\lim_{\alpha\to\infty}\sup_{k\ge1}\mbf{P}\{\<\psi_k, Y_t^k\>\ge
\alpha\} = 0.
 \eeqnn
Then a criterion of Ref.~\cite{A78} implies that the sequence
$\{\<\psi_k,Y_t^k\>: t\ge 0\}$ is tight in $C([0,\infty), \mbb{R})$.
Let $\Phi_p$ be defined by (\ref{5.0}). By Theorem~3.7.1 of
Ref.~\cite{D93}, $\{k^{-1}\Phi_p Y_t^k: t\ge 0\}$ is a tight
sequence in $C([0,\infty), M(\bar\mbb{R}))$. Now the result of the
first step implies that $\{k^{-1}\Phi_pY_t^k: t\ge 0\}$ converges to
$\{t\Phi_p\lambda: t\ge 0\}$ in probability in $C([0,\infty),
M(\bar\mbb{R}))$. Since all the processes live in $C([0,\infty),
M(\mbb{R}))$, the convergence also holds in this smaller space. In
other words, $\{k^{-1}Y_t^k: t\ge 0\}$ converges to
$\{q(\infty)t\lambda: t\ge 0\}$ in probability in $C([0,\infty),
M_p(\mbb{R}))$.  \qed

The above theorem implies the following scaling limit theorem for
the local time of the ISDSM.

\bcorollary\label{c5.1} Let $z(\cdot,\cdot)$ denote the local time
of $\{Y_t: t\ge0\}$ given by the second term on the right hand side
of (\ref{4.10}). Under the condition of Theorem~\ref{t5.1},
$z_k(t,\cdot) := k^{-5}z(k\,\cdot,k^2t)$ converges weakly to
$q(\infty)t^2/2$ in probability as $k\to \infty$. \ecorollary

For a finite reference immigration measure, we can prove a limit
theorem which gives an interesting random limit process. To describe
the limit process we introduce the following concept: A
two-parameter process $\{y_r(t): t\ge r\ge0\}$ is called a
\textit{restricted coalescing Brownian flow} (RCBM flow) with speed
$\rho>0$ provided
 \begin{itemize}

\item[(5.a)]
for any fixed $r\ge0$, the process $\{y_r(t): t\ge r\}$ is a
Brownian motion with speed $\rho$ started from $y_r(r)=0$;

\item[(5.b)]
for any fixed $s\ge r\ge 0$, the process $\{y_s(t) - y_r(t): t\ge
s\}$ is a Brownian motion with speed $2\rho$ stopped at zero.
 \end{itemize}
Let $m\in M(\mbb{R})$ and let $q(\cdot,\cdot)$ be a bounded Borel
function on $M(\mbb{R})\times \mbb{R}$ satisfying the local
Lipschitz condition (\ref{3.3}) for $p=0$. In addition, we assume
 \begin{itemize}

\item[(5.c)]
there is a constant $\epsilon>0$ so that $q(\nu,a) \ge \epsilon$
for all $\nu \in M(\mbb{R})$ and $a \in \mbb{R}$;

\item[(5.d)]
$q(\nu,a) \to q(a)$ as $\<1,\nu\> \to\infty$ for all $a \in
\mbb{R}$.

 \end{itemize}
Suppose that $\{W(ds,dy)\}$ and $\{N_1(ds,da$, $du,dw)\}$ are given
as in (2.a) and (2.d) and let $\{Y_t: t\ge0\}$ be defined by
(\ref{5.1}). Let $Y_t^k(dx) = k^{-2}Y_{k^2t}(kdx)$.

\btheorem\label{t5.2} Under the above conditions, as $k\to
\infty$, the sequence $\{Y_t^k: t\ge0\}$ converges in distribution
on $C([0,\infty), M(\mbb{R}))$ to
 \beqlb\label{5.6}
Y_t^\infty
 :=
\int_0^t\int_{\mbb{R}}\int_0^{q(a)}\int_{W_0} w(t-s)
\delta_{x^\infty_s(t)} N_1(ds,da,du,dw), \qquad t\ge0,
 \eeqlb
where $\{x^\infty_s(t): t\ge s\}$ is an RCBM flow independent of
$\{N_1(ds,da,du,dw)\}$ with speed $\rho(0)$. \etheorem

\proof Let $\{x_{s,b}^k(t): t\ge s\}$ be defined as in the proof
of Theorem~\ref{t5.1}. For any $\phi\in C^2(\mbb{R})$ we have
 \beqlb\label{5.7}
\<\phi,Y_t^k\>
 =
\int_0^t\int_{\mbb{R}}\int_0^{q(Y_{k^2s},a)}\int_{W_0} w_k(t-s)
\phi(x^k_{s,a/k}(t)) N_1(k^2ds,da,du,dw).
 \eeqlb
It is easy to see that $N_1(k^2ds,da,du,dw)$ is a Poisson random
measure with intensity $k^2dsm(da)du \mbf{Q}_\kappa(dw)$. Let
$X_t^k$ be defined by the right hand side of (\ref{5.1}) with
$q(Y_s,a)$ replaced by $q(a)$. Then we have
 \beqlb\label{5.8}
\<\phi,X_t^k\>
 =
\int_0^t\int_{\mbb{R}}\int_0^{q(a)}\int_{W_0} w_k(t-s)
\phi(x^k_{s,a/k}(t)) N_1(k^2ds,da,du,dw).
 \eeqlb
The scaling property of the $\sigma$-branching diffusion implies
that
 \beqnn
\<\phi,X_t^k\>
 =_d
\int_0^t\int_{\mbb{R}}\int_0^{q(a)}\int_{W_0} w(t-s)
\phi(x^k_{s,a/k}(t)) N_1(ds,da,du,dw).
 \eeqnn
Therefore the distribution of $\{\<1,X_t^k\>: t\ge 0\}$ on
$C([0,\infty), \mbb{R}_+)$ is independent of $k\ge1$. Indeed, by
Theorem~4.1 of Ref.~\cite{PY82} it is easy to see that
$\{\<1,X_t^k\>: t\ge 0\}$ is distributed as the $\sigma$-branching
diffusion with immigration that solves the stochastic equation
 \beqnn
dz_t
 =
\sqrt{\sigma z_t}dB_t + \<q,m\>dt, \qquad z_0=0.
 \eeqnn
Given $(s_1,\cdots,s_m; a_1,\cdots,a_m)$, it can be proved as
Theorem~2.2 of Ref.~\cite{DLZ04} that $\{(x_{s_1,a_1/k}^k (\cdot),$
$\cdots, x_{s_m,a_m/k}^k(\cdot))\}$ converges in distribution to
$\{(x^\infty_{s_1}(\cdot), \cdots, x^\infty_{s_m}(\cdot))\}$. By
dominated convergence it is easy to show that $\{\<\phi,X_t^k\>:
t\ge 0\}$ converges in finite dimensional distributions to
 \beqnn
\<\phi,Y_t^\infty\>
 =
\int_0^t\int_{\mbb{R}}\int_0^{q(a)}\int_{W_k} w(t-s)
\phi(x^\infty_s(t)) N_1(ds,da,du,dw), \qquad t\ge0.
 \eeqnn
On the other hand, from (\ref{5.7}) and (\ref{5.8}) it is easy to
see that
 \beqnn
\big|\<\phi,X_t^k\> - \<\phi,Y_t^k\>\big|
 \le
\int_0^t\int_{\mbb{R}}\int_{q(a)\land q(Y_{k^2s},a)} ^{q(a)\vee
q(Y_{k^2s},a)} \int_{W_0} w_k(t-s) \|\phi\| N_1(k^2ds,da,du,dw)
 \eeqnn
and so
 \beqlb\label{5.9}
\mbf{E}\big[\big|\<\phi,X_t^k\> - \<\phi,Y_t^k\>\big|\big]
 \le
\|\phi\|\int_0^tds\int_{\mbb{R}} \mbf{E}[|q(a) - q(Y_{k^2s}, a)|]
m(da),
 \eeqlb
By Condition (5.c) we find that $\<1, Y_{k^2s}\> = k^2 \<1,
Y_s^k\> \to \infty$ in probability for every $s>0$. Then
(\ref{5.9}) and Condition (5.d) implies that $|\<\phi,X_t^k\> -
\<\phi,Y_t^k\>| \to 0$ in probability. Therefore $\{\<\phi,
Y^k_t\>: t\ge 0\}$ also converges to $\{\<\phi, Y^\infty_t\>: t\ge
0\}$ in finite dimensional distributions. The tightness of
$\{Y^k_t: t\ge 0\}$ can be established as in the proof of
Theorem~\ref{t5.1}, so the sequence converges to $\{Y^\infty_t:
t\ge 0\}$ in distribution on $C([0,\infty), M(\mbb{R}))$. \qed

\bcorollary\label{c5.2} Let $z(\cdot,\cdot)$ denote the local time
of $\{Y_t: t\ge0\}$. Under the condition of Theorem~\ref{t5.2},
$k^{-4}z(k\,\cdot,k^2t)$ converges weakly in distribution to
 \beqlb\label{5.10}
z^\infty(\cdot,t)
 :=
\int_0^t\int_{\mbb{R}}\int_0^{q(a)}\int_{W_0}N_1(ds,da,du,dw)
\int_s^t w(v-s)l_s^\infty(\cdot,dv),
 \eeqlb
where $\{l_s^\infty(b,u): u\ge s\}$ denotes the local time of
$\{x^\infty_s(u): u\ge s\}$ at $b\in \mbb{R}$. Moreover, the
process $\{z^\infty(b,t): b\in \mbb{R}, t\ge 0\}$ has a version
which is H\"older continuous of order $\alpha$ for every $\alpha <
1/2$. \ecorollary

\proof The first part is an immediate consequence of
Theorem~\ref{t5.2}. The second part follows by arguments similar
to those given in the proof of Theorem~\ref{t4.1}. \qed

\bigskip\bigskip

\centerline{\large\bf Acknowledgements}

\bigskip

Z.L. was supported by NSFC (10121101 and 10525103) and NCET
(04-0150). J.X. was supported partially by NSA.

\end{document}